\RequirePackage{fix-cm}
\documentclass[smallcondensed]{svjour3}     
\smartqed  
\usepackage[a4paper]{geometry}
\usepackage{graphicx}
\usepackage{comment}
\usepackage{url}
\usepackage{amsmath,amssymb,amsfonts}
\usepackage{algorithmicx}
\usepackage{algpseudocode}
\usepackage{algorithm}
\algnewcommand{\IIf}[1]{\State\algorithmicif\ #1\ \algorithmicthen}
\algnewcommand{\EndIIf}{\unskip\ \algorithmicend\ \algorithmicif}
\makeatletter
\algnewcommand{\LineComment}[1]{\Statex \hskip\ALG@thistlm #1}
\makeatother
\usepackage{hyperref}
\usepackage{mathtools}
\usepackage[misc,geometry]{ifsym}
\usepackage[]{natbib}
\usepackage{mathptmx}
\usepackage{graphicx}
\usepackage{caption}
\usepackage{subcaption}
\usepackage{geometry}
\usepackage{pdflscape}
\usepackage{booktabs}
\usepackage{epstopdf}
\usepackage{tabularx}
\usepackage[table]{xcolor}
\graphicspath{{./Figures/}}

\usepackage{color}

\begin{document}

\title{A hybrid swarm-based algorithm for single-objective optimization problems involving high-cost analyses} 

\titlerunning{A hybrid swarm-based algorithm for optimization problems involving high-cost analyses} 

\author{Enrico Ampellio \and Luca Vassio}

\institute{E. Ampellio (\Letter) \at 
              Dipartimento di Ingegneria Meccanica e Aerospaziale, Politecnico di Torino, 10129 – Torino (TO), Italy \\
              \email{enrico.ampellio@polito.it} 
           \and
           L. Vassio \at 
           Dipartimento di Elettronica e Telecomunicazioni, Politecnico di Torino, 10129 – Torino (TO), Italy \\
              \email{luca.vassio@polito.it}
}

\date{Received: May 26, 2015 / Accepted: date}

\maketitle

\begin{abstract}

In many technical fields, single-objective optimization procedures in continuous domains involve expensive numerical simulations.
In this context, an improvement of the \textit{Artificial Bee Colony} (ABC) algorithm, called the \textit{Artificial super-Bee enhanced Colony} (AsBeC), is presented.
AsBeC is designed to provide fast convergence speed, high solution accuracy and robust performance over a wide range of problems. 
It implements enhancements of the ABC structure and hybridizations with interpolation strategies. The latter are inspired by the quadratic trust region approach for local investigation and by an efficient global optimizer for separable problems. Each modification and their combined effects are studied with appropriate metrics on a numerical benchmark, which is also used for comparing AsBeC with some effective ABC variants and other derivative-free algorithms.
In addition, the presented algorithm is validated on two recent benchmarks adopted for competitions in international conferences. Results shows remarkable competitiveness and robustness for AsBeC.

\keywords{modified Artificial Bee Colony \and engineering optimization \and interpolation strategies \and algorithm comparison}

\end{abstract}

\section{Introduction}\label{Intro}

Numerical optimization in applied sciences usually concerns simulation-based problems in which computational analyses are expensive in terms of time and resources (S. Rao, \citeyear{RA09}; K. Deb, \citeyear{DE12}). In engineering, complex non-linear models are commonly used (Lagaros and Papadrakakis, \citeyear{LP15}), e.g., computational fluid dynamics and finite element methods. This field demands fast algorithms for solving optimization problems using as few function evaluations (FEs) as possible. The present work focuses on single-objective problems with low, but not minimal dimensionality. 

In the literature, there are several architectures and algorithms suitable for global optimization in tough numerical problems (Floudas and Pardolos, \citeyear{FP09}). Among the most studied are surrogates (also called meta-models), heuristic and statistical methods, trust region approaches and nature-based algorithms.

Surrogate models (Koziel et al., \citeyear{KO13}) are effectively exploited in a large class of real-world applications. They are often the first choice and sometimes the only practical one when very expensive analyses are involved. Response surfaces (Box and Draper, \citeyear{BD07}), artificial neural networks (Iliadis and Jayne, \citeyear{IJ15}), radial basis functions (Buhmann, \citeyear{BU03}) and Kriging (Simpson et al., \citeyear{SI01}) are surrogates widely adopted to support high computing load optimizations. One of the first optimizers fully based on meta-models is the Efficient Global Optimization (Jones et al., \citeyear{JO98}). This algorithm and its successors are expressively targeted to problems where function evaluations are strictly limited by time or cost. Heuristic and statistical techniques such as Sequential Parameter Optimization (SPO, Bartz-Beielstein et al., \citeyear{BA05}) are also related to meta-models and intensively used in applicative fields. They are focused on selecting and tuning both the optimization methodology and the search methods, but they imply that efficient algorithms and/or surrogates for a specific problem are a priori available. Nonetheless, ensuring the reliability of surrogates is still a sensitive issue and the training of meta-models could not be cost-effective when considering highly non-linear, noisy, multi-modal and discontinuous problems. Moreover, this approach becomes computationally impracticable when increasing the problem dimensionality.
As a result, it is attractive to develop direct search algorithms designed for costly optimizations, especially to cope with higher dimensions (e.g., $10$ dimensions).

Among direct search algorithms, derivative-free trust region optimizers (Conn et al., \citeyear{CO00}), such as the quadratic UOBYQA (Powell, \citeyear{PO00}) and related approaches (e.g., CONDOR by Berghen and Bersini, \citeyear{BB05}; BOBYQA by Powell, \citeyear{PO09}), are effective for heavy computing loads and noisy objective functions. 
The hybridization with quadratic local meta-models is also used to tackle costly problems, for instance in evolutionary strategies (lmm-ES by Kern et al., \citeyear{KE06}), in which surrogates are periodically called to improve performance and save expensive function evaluations (FEs).

Meta-heuristic direct-search algorithms play a predominant role when dealing with single-objective optimizations that comprise many continuous variables and noisy multi-modal target functions, and where analytical properties such as differentiability are not known a priori (Talbi, \citeyear{TA09}). Several nature-based techniques have been developed such as genetic algorithms (Goldberg, \citeyear{GO89}), particle swarm optimization (Kennedy and Eberhart, \citeyear{KE95}), ant colony optimization (Dorigo and St{\"u}tzle, \citeyear{DO04}), differential evolution (DE by Price et al., \citeyear{PR05}) and evolutionary strategies (Beyer and Schwefel, \citeyear{BS02}).

One of the more recent and promising nature-inspired meta-heuristic methods is the Artificial Bee Colony algorithm (ABC) by Karaboga and Akay (\citeyear{KA07}). ABC has been compared with other nature-inspired algorithms and it showed good results for a great variety of optimization problems (Karaboga and Akay, \citeyear{KA09}) as well as for real-world applications (Akay and Karaboga, \citeyear{AK12}). Besides, a lot of test case validations and advancements have been done in recent years (Bolaji et al., \citeyear{BO13}). The present paper further improves the accuracy and convergence speed of ABC algorithms to make them suitable for continuous optimization problems with expensive analyses. This means obtaining better solutions using the same maximum number of function evaluations. The main idea behind the proposed algorithm is to hybridize ABC with the basic principles of derivative-free, quadratic trust region optimizers for high-cost objective functions. Moreover, another hybridization with the fast and robust  S.T.E.P method (Swarzberg~et~al.~\citeyear{SW94}) is introduced to efficiently cope with separable functions. In the end, AsBeC aims at showing robust performance over a wide range of problems by successfully merging the main advantages of swarm-based and interpolative methods.

To introduce the main peculiarities of ABC, a brief review of the original algorithm and its popular modifications is presented in Section~\ref{ABC}. The proposed enhancements and hybridizations leading to the novel Artificial super-Bee enhanced Colony algorithm (AsBeC) are then discussed in Section~\ref{Tech}. 
The analytical benchmarks and new metrics used for comparisons are respectively introduced in Section~\ref{benchmark} and Section~\ref{metrics}. The performance of the new algorithm is accurately evaluated and compared to  ABC through extensive simulations in Section~\ref{Test}. In Section~\ref{Comparisons}, some ABC variants and other state-of-the-art algorithms for direct search and appropriate for costly optimizations are considered for comparisons. In Section~\ref{Validation}, the algorithm is validated using two other famous benchmarks for expensive optimization. Conclusions and insights about future work on the topic are presented in Section~\ref{Conclusion}. An electronic appendix gathers all comparative results in expanded form.

\vspace{-5pt}
\section{ABC and its modifications in the literature}\label{ABC}

The Artificial Bee Colony (ABC) algorithm by Karaboga (\citeyear{KA07}) is an example of swarm intelligence as it models the collective behaviour of decentralized and self-organized systems consisting of simple agents. The algorithm reproduces the behaviour of a honey bee colony searching the best nectar source in a target area. Some bees, called the employees, are each assigned to a nectar source and search the space near it. Then, they come back to the hive and communicate to other bees, called the onlookers, the position of the best nectar sources found through a dance (Frish, \citeyear{FR67}), so that onlookers can help the employed bees in searching the most promising regions. Nectar sources that reveal themselves non-productive are abandoned in favour of new positions randomly chosen in the search domain, which are investigated by scout bees.

In the optimization context, nectar sources represent the seed solutions and the nectar amount corresponds to the objective function value to optimize; "non-productive" sources are those not improving for some time, controlled by the limit parameter. An entire optimization sequence composed by the employee group, performing the exploration, the onlooker group, performing the exploitation, and the scout bee, maintaining the population diversity, is called a \textit{cycle}. For an ABC pseudo-code refer to Section~2.1 in Karaboga and Akay (\citeyear{KA07b}).

The basic search movement for a seed solution $x^j$ is called mutation equation and it modifies only one variable at a time, $i$, chosen uniformly at random. Another seed solution $x^k\neq x^j$ is also chosen uniformly at random and the new candidate solution $x^{j,rnd}$ in variable $i$ is:

\begin{equation}\label{eq.step}
	x^{j,rnd}_i=x^j_i+rnd[-1,1]\cdot(  x^j_i-x^k_i	)
\end{equation}

where $rnd[-1,1]$ is a uniformly distributed random number in the continuous range $[-1,1]$.

The ABC algorithm tries to balance exploration and exploitation behaviours, offering worthy global and local search skills at once. As a result, when compared with other methods, ABC demonstrates high-performance and suitability for a broad spectrum of optimization problems. To enter in details consult Karaboga (\citeyear{KA07}) and for extensive simulation comparisons refer to Karaboga and Akay (\citeyear{KA09}).
According to many authors, ABC algorithm is simple to implement, easy to be effectively parallelized (e.g., Subotic et al., \citeyear{SU11}) and hybridized, driven by few control parameters, flexible and robust (e.g., Karaboga et al. \citeyear{KA14}). On the other hand, its semi-random movement is not taking into account the local shape of the function. Moreover, the algorithm does not exploit the history of the tested solutions.

Since the original paper by Karaboga (\citeyear{KA07}), many researches were developed to overcome the deficiencies of the original ABC. In particular, many scientists focused on the common goal to advance the local search mechanism without worsening the global one. In this regard, specific modifications of various nature have proven to be convenient (Bolaji et al., \citeyear{BO13}).
Among them, the popular ABC variants by Gao and Liu (e.g., IABC, \citeyear{GL11} and BABC, \citeyear{GL12}) and by Kong et al. (e.g., \citeyear{KON13}) are focused on initialization enhancements and DE/PSO inspired mutation equations. Nevertheless, many preliminary tests performed by the authors of the present paper showed that these modifications are not of major impact when dealing with a limited number of FEs and small colonies. Their performance in this framework is comparable with the simpler and well-known GABC (Zhu and Kwong, \citeyear{ZK10}), which drives the classic bee movement towards the global best solution. 

Some other variants give more importance to the exploitation attitude by altering the fitness function or the assignation methodology for onlookers (e.g., Subotic, \citeyear{SU11b}). Among them, JA-ABC5 (Sulaiman et al., \citeyear{SL15}) introduces a multi-stage re-organization of the bees which seems particularly beneficial. 

Many modifications have been proposed coupling ABC either with classical optimization methods for local search such as simplex (HSABC by Kang et al., \citeyear{KAN09}), the Rosenbrock method (RABC  by Kang et al., \citeyear{KAN11}) and pattern search (HJABC by Kang et al., \citeyear{KAN13}), or with other robust and efficient well-known techniques for global search, such as genetic algorithms (GA-ABC by Zhao et al., \citeyear{ZH10}) or evolutionary strategies (BE-ABC by Li and Li, \citeyear{LL12}). Among these methods, the RABC based on the Rosenbrock directional search (Rosenbrock, \citeyear{RO60}) shows fast convergence to the global optimum.

Despite the great amount of available literature on ABC variants, to the best of our knowledge no paper underlines a performance gain even with very few function evaluations. As well, no hybridization with methods similar to trust region techniques has been analysed so far, even though they are recognized to be among the fastest techniques, especially when dimensionality is small (Rios and Sahinidis, \citeyear{RS13}). These reasons justify the introduction and analysis of ABC modifications specifically targeted to costly optimization and real-world oriented problems.

\vspace{-5pt}
\section{Improving techniques for the AsBeC algorithm}\label{Tech}

Among all the proposed modifications of the ABC algorithm in the literature, some techniques are on average more effective or more robust (Liao et al.,\citeyear{LA13} and Aydin, \citeyear{AY15}) for specific problems.
In the present work, the authors discover that modifying the onlooker assignment stage, exploiting the data history according to the theory of collective memory in social animals (Couzin et al., \citeyear{CO02}) and applying some hybridizations like local interpolation are especially beneficial when the number of FEs is very limited. The presented techniques have the purpose of achieving the fastest solution improvement, without clustering the swarm and leading to premature convergence at the same time. 

The sections below describe each technique in details. A minimization problem is considered without loss of generality. The final AsBeC pseudo-code is outlined and commented in Section~\ref{AsBeC}.

\vspace{-5pt}
\subsection{Strictly biased onlooker assignment}\label{Tech1}

In ABC, onlookers are assigned to seed solutions by a stochastic rule. A certain probability is associated to each of them: it is proportional to the solution quality $f$, through a fitness function \emph{fit}. Probability and fitness definitions are taken from Akay and Karaboga (\citeyear{AK10}). However, if the value of $f$ for all the seed solutions is similar, that fitness function returns almost identical values and the relative probabilities are alike. Hence, the standard fitness formulation is not able to distinguish the solution quality in relative terms. In order to always strengthen exploitation, it is possible to unbalance onlooker assignation by normalizing the fitness \emph{fit} to a new re-scaled function ${\emph{fit}}^{res}$:

\begin{equation*}
{\emph{fit}}^{res}(x^k)=\frac{\emph{fit}(x^k)-\min_{s=1,\ldots,SN}⁡{fit(x^s)}}{\max_{s=1,\ldots,SN}⁡{\emph{fit}(x^s)}-\min_{s=1,\ldots,SN}⁡{\emph{fit}(x^s)}}
\quad \forall k \in \{1,\ldots,SN\}
\end{equation*}

where $x^k$ is a seed solution and $SN$ is the number of seed solutions. The above formulation for ${\emph{fit}}^{res}(x^k)$ is generalized and still valid even when changing the definition of \emph{fit}. Different from the rank-based fitness transformation of Kang et al. (\citeyear{KAN11} and \citeyear{KAN13}), it takes into account also the quality difference among seed solutions. In fact, $\emph{fit}^{res}$ scores the seed solutions from 0 to 1 according to their distance from the best.

In the original ABC formulation, there is a possibility to assign few or no bees to the best seed solutions. The previous scenario is unlikely when dealing with a large-size swarm; however, for small-size colonies a different approach should be preferred in order to save function evaluations. A deterministic rule for onlookers assignation is proposed in place of the original stochastic one, in order to associate a prearranged number of bees $n^k$ to seed solution $k$ proportionally to its ${\emph{fit}}^{res}$:

\begin{equation*}
n^k=\Bigg{\lfloor} \frac{ON\cdot{\emph{fit}}^{res}(x^k)}{\sum_{s=1}^{SN}{\emph{fit}}^{res}(x^s)}\Bigg{\rfloor} \quad \forall k \in \{1,\ldots,SN\}
\end{equation*}

where $ON$ is the onlooker number and the operator $\lfloor \cdot \rfloor$ returns the largest integer smaller than the argument. In case the number of onlookers assigned, $\sum_{k=1}^{SN} n^k$, is lower than the total number of onlookers $ON$, the difference $ON-\sum_{k=1}^{SN} n^k$ is assigned to the best seed solution. In practice, no bees are ever sent to the worst seed solution and at least one is always sent to the best one. The application of this deterministic rule is similar to the idea of exploiting only few best solutions of the JA-ABC5 algorithm and it always encourages the repositioning near the global best, as in GABC; however, it is self-adapted during the optimization process.

Notice that this strictly biased onlooker (BO) assignment, through fitness re-scaling and deterministic rule, does not appreciably affect the global search skills of the swarm. In fact, exploration is essentially up to the employee movement, which remains unaltered. This technique intentionally disregards the worst areas, similarly to the principle driving other swarm-based algorithms such as the FIREFLY algorithm (Yang, \citeyear{YA09}). 

\vspace{-5pt}
\subsection{Postponed onlooker dance}\label{Tech2}

This modification enables the onlooker bees to seek nectar near the best regions for a longer time. The entire onlooker group performs more than one iteration sequentially, without altering the updating order of seed solutions inside each iteration. As in the original ABC, a seed solution is immediately updated whenever a bee finds an improvement for it. In the postponed onlooker dance (PD) iterations, seed solutions are never reassigned, but their updated versions are the starting configurations for the new movements. In this way, each onlooker bee has multiple chances to improve the seed solution to which it was assigned and to interact with the other proposed techniques. 
Anyway, the number of PD iterations should be maintained small in order not to affect exploration and decelerate the convergence speed to the global optimum. In this paper, the number of reiterations on the onlooker group is set to three.

This modification acts like enlarging the onlooker group, but repeating their movements instead of increasing their number. Hence, the best seed solutions have a higher probability to improve and to be selected for mutation than the others, as onlookers are not distributed on a larger basis.

This postponed onlooker dance among bee groups is inspired by nature. In fact, bees normally come back to the hive to communicate their seed solutions only after some time, when there is more probability to have collected new pollen. This technique is inspired by the work of Subotic (MO-ABC, \citeyear{SU11b}), except that the PD proposed is a deterministic procedure.  

\vspace{-5pt}
\subsection{Local interpolation}\label{Tech3}

In order to quickly improve the best solutions in their local neighbourhood, the concepts of opposition principle and parabolic interpolation are introduced.

The opposition-based learning for the ABC algorithm has already been suggested in the OABC by El-Abd (\citeyear{EL11}). In the present paper, the following version is proposed: whenever a pseudo-random movement for the seed solution $k$ to a new position $x^{k,rnd}$ does not produce any improvement, it is possible to move in the same direction with the same step, but in the opposite sense. If this new position $x^{k,opp}$ is better than the previous, the correspondent seed solution $k$ is updated. As a result, the opposite principle is a linear local estimator when the step $\|{x^k-x^{k,rnd}}\|$ is small enough and the objective function is continuously differentiable in the neighborhood of $x^k$. The operator $\|{ \cdot }\|$ expresses the Euclidean norm.

On the other hand, the parabolic interpolation estimates the local and directional curvature of the objective function, acting as a second order optimization method with partial Hessian computation. This technique follows the opposition principle: if the opposite step is not successful, the multi-dimensional parabola passing through three already known points, i.e., (i) the point $x^k$, (ii) the point generated by the first random movement $x^{k,rnd}$ and (iii) the opposite point $x^{k,opp}$, is calculated and its minimum is tested. The opposition principle and second order interpolation are performed in excluding sequence: if the first random movement is improving, none of the two following steps is carried out; if the opposite movement is improving, the parabola will not be estimated; if the two previous steps lead to solutions worse than the seed solution, then the second order interpolation is computed.

Local interpolation (LI) is used in ABC when at least three onlooker bees are assigned to the same seed solution. When the PD technique is activated, the onlooker group has more chances to exploit local interpolation.
During the onlooker phase, the previous bee movements relative to the same seed solution are recorded and shared among bees, also considering different PD iterations.
When performing local interpolation, it is possible to use candidate solutions generated by different onlooker bees that are assigned to the same seed solution.
Using three iterations for the PD technique, a seed solution with just one onlooker assigned could still perform a complete sequence of random, opposite and parabolic interpolation movements.

The raw approximation offered by local interpolation is balanced by the low number of function evaluations required, i.e., respectively one and two additional FEs for the opposite movement and the parabolic interpolation, regardless of the problem dimensionality.

\vspace{-5pt}
\subsection{Quadratic prophet}\label{Tech4}

A hybridization with trust region concepts is introduced in order to improve the convergence speed, especially in the early phases. The quadratic prophet (QP) hybridizes ABC with a fast and robust local search method, which is the same philosophy of other hybrids such as RABC (Kang et al., \citeyear{KAN11}).
The key idea is to exploit the function evaluations already performed to create multiple quadratic polynomial interpolators around the best regions. Hence, these quadratic models are built in the neighborhood of the seed solutions by using their closest samplings in the so-called \textit{hive memory}. This memory consists in all the solutions generated by the AsBeC algorithm so far. The global minima of the quadratic models, if they exist, are then tested. 
The lmm-CMAES algorithm by Kern et al. (\citeyear{KE06}) is similar, but it uses the quadratic models as response surfaces to be called in place of direct function evaluations. 
Instead, QP does not require additional FEs to populate the samplings for the interpolations; its final cost is at the most equal to the number of seed solutions if all quadratic models succeeded. Whenever the QP finds a better solution, the related seed solution is immediately updated.

Although this method does not regulate directly the trust region radius, it adapts itself in accordance to the hive memory. In fact, if the samplings around a given seed solution are dense, the algorithm is exploiting the area, and thus the trust region radius is small in order to refine the local solution. On the contrary, when the samplings are sparse, the trust region radius is large to support exploration. This self-adaptation resembles the mechanism used in the CMA-ES algorithm  (Hansen, \citeyear{HA06}).
Other famous trust region methods, such as the classic Levenberg-Marquardt (Levenberg, \citeyear{LE44}; Marquardt, \citeyear{MA63}) or the recent BOBYQA (Powell, \citeyear{PO09}) algorithms follow the same basic principle. The main difference is that the QP technique uses the data that are already available within the AsBeC algorithm, instead of adopting a costly and complex population strategy to train and use quadratic approximations. On the other hand, the quality of the QP model depends on the quality of the samplings, hence on the ability of AsBeC to properly investigate the search space surrounding the minima. Therefore, the Quadratic Prophet technique can hybridize any optimization algorithm able to provide well-distributed samplings.

The implementation of QP is based on quadratic polynomial surfaces, whose complete form, $Q_C$, and reduced form, $Q_R$, for a generic problem in $D$ dimensions are expressed as:

\begin{equation*}
Q_C(x_1,x_2,...,x_D)=\sum_{i=1}^{D}\sum_{j=1}^{D}a_{ij}\cdot x_i \cdot x_j + \sum_{i=1}^{D} b_i \cdot x_i + R \qquad 
Q_R(x_1,x_2,...,x_D)=\sum_{i=1}^{D}a_{ii}\cdot x_i^2 + \sum_{i=1}^{D} b_i \cdot x_i + R
\end{equation*}

where $a_{ij}$, $b_i$ and $R$ are the $C$ coefficients of the model, while $x_i$ are the independent variables of the problem. $Q_C$ needs  $C=0.5\cdot (D+2) \cdot (D+1)$ samplings to be exactly defined, since $a_{ij}=a_{ji}$. The reduced quadratic surface $Q_R$ , without any mixed term, only needs $2 \cdot D+1$ samplings.
The reduced model is not able to capture rotated functions, but it is enough to cope with bowl shape regions. In order to activate this technique as soon as possible, reduced models $Q_R$ are built starting when the algorithm provides at least  $2 \cdot D+1$ solutions, ending up with complete models $Q_C$ when more than $C$ samplings are available.

The model coefficients are obtained in the simplest possible way, that is by solving exactly the linear system of equations for the $C$ samplings nearest to a given seed solution. If the Hessian matrix of the model is diagonal positive, then the point in which the Jacobian determinant is zero is tested. This point could represent either a local minimum, if the Hessian is positive definite, or a saddle point dominated by squared terms. Actually, the latter can be wrongly recognized as a saddle instead of a minimum due to inaccuracies. Nonetheless, it represents a key point to better describe the objective function contour in that region. Anyway, eventual redundant QP predictions are not tested twice, in order to save precious FEs.

Notice that the model could misinterpret the local behaviour of the objective function if the sampling is ill-conditioned or if $f$ is not quadratic or noisy in that area, as common in real-world problems. Then, if the exact solution does not lead to a point to test, the system is solved by means of an iterative method with large tolerance on the same $C$ samplings. This alternative is aimed at obtaining quickly a simple-to-handle coarse shape interpreter, passing nearby the sample points and not exactly through them. This inaccurate model could still be able to capture the objective function profile and drive the search towards the most promising regions. The tolerance for the iterative resolution of quadratic systems is set to $10^{-1}$ by using the QMR method (Roland and Nachtigal, \citeyear{RN91}).

In order to increase the likelihood of a well-conditioned sampling, the points should be heterogeneously different in their coordinates. There are two main features that affect the sampling diversity: the standard bee movement, changing only one variable at a time, and the LI technique, producing sequentially aligned points. Therefore, when the QP is activated, LI is deactivated and the mutation equation is slightly modified.
Following the idea introduced by Akay and Karaboga (\citeyear{AK10}) for the MABC, more than one variable is chosen when moving a bee, according to the standard ABC equation.
The number of variables to be changed concurrently is set as half of the total number of dimensions, i.e., $D/2$.
Clearly, these two modifications overweight exploration and hence have to be restricted to the early phases of the optimization. In this work, they are adopted while the total number of solutions tested is lower than $2\cdot C= (D+1) \cdot (D+2)$. In this way, a large pool of samplings more likely to be well-conditioned is available for complete quadratic models.

This QP technique is a finer exploitation of the seed solutions and, hence, it is implemented as an additional exploitative stage at every cycle, after the onlooker phase. It is also used once during the initialization. Indeed, the number of random initial evaluations is set to be (at least) equal to $2 \cdot D+1$, and a fully reduced $Q_R$ is already feasible. 

It must be clear that the QP is designed to deal with problems in which the objective function evaluation represents the bottleneck of computational times. In fact, the computational cost of this technique is dominated by the resolution of the linear system, that in its general form requires $O(C^3)$ operations at each cycle, for each food source. Given that $C$ increases along with $D$, the technique becomes numerically heavy when addressing high dimensionality problems. 

\vspace{-5pt}
\subsection{Systematic global optimization}\label{Tech5}

The S.T.E.P algorithm by Swarzberg~et~al.~(\citeyear{SW94}) shows high convergence rate for highly multi-modal and complex functions of one variable. It is derivative-free, not population-based, without parameters to tune and it does not assume any property for the function to optimize. It relies on the history of the points tested during its functioning. 
The domain boundaries and a random point between them, called the \textit{context}, are evaluated. These three points define two \textit{partitions}. For each partition is computed the curvature a parabola should have to enclose the best so far solution inside it. This curvature is called the partition \textit{difficulty}. The S.T.E.P method iterates by sampling the center of the partition with smallest difficulty, i.e. the one identified with the greatest chance of improving the best solution found so far, until a given tolerance is reached. The S.T.E.P. authors showed that around 50 FEs are sufficient to solve difficult univariate problems, such as Michalewicz's second function.
The simplest multivariate extension to the S.T.E.P. algorithm solves one dimension at a time. The recent version by Baudi\v{s} and Po\v{s}\'{\i}k~(\citeyear{PO15}) interleaves the steps of the univariate solvers such that all dimensions are optimized concurrently, in a similar way to Rosenbrock algorithm (\citeyear{RO60}). Each time, the dimension to investigate is chosen following a round-robin scheduling.

Given the fast convergence and robust performance of the S.T.E.P. algorithm, it can be hybridized with ABC, which is more oriented to exploration. However, in a multivariate environment S.T.E.P. based solvers are very effective only on additively separable or quasi-separable problems, in which the correlations among variables are weak or only few variables are correlated. As a consequence, this technique is activated if a test for quasi-separability is successfully passed. The proposed test investigates the separability in the cheapest possible way, giving just a necessary but not sufficient condition. 
For an additively separable function, each first-order partial derivative ${\partial f}/{\partial x_i}$ is independent from other variables $x_j \neq x_i$, giving ${\partial^2 f}/{\partial x_i \partial x_j=0}$. For each design variable $x_i$, another variable $x_j$ with $i \neq j$ is randomly chosen to asses the correlation between them.
Then, the partial derivative ${\partial f}/{\partial x_i}$ is estimated for two different values of $x_j$,  while all other coordinates are kept constant.
If the two estimations vary below a given tolerance, here chosen as $10^{-3}$, the two variables could be quasi-separable. The same check is repeated for ${\partial f}/{\partial x_j}$ for two different values of $x_i$. In order to estimate the aforementioned derivatives through variations, four FEs are needed. These four solutions are chosen around the active seed solutions.  
If the non-separability condition on the derivatives holds for D couples of variables, the test is passed. The procedure is interrupted whenever a variable is found to be correlated. This test needs at most $2 \cdot D + 1$ function evaluations, since a complete set of partial derivatives in the same point can be estimated only once and reused for all the couples.

If this test is successfully passed,  the interleaved solver by Baudi\v{s} and Po\v{s}\'{\i}k~(\citeyear{PO15}) is run as a one-shot technique after the initialization phase, using the best so-far solution as context and investigating all the variables, in deterministic sequence. For this reason, the proposed hybridization is called systematic global optimization (SGO). At the end of the SGO phase, the best seed solution is moved to the best solution found. The maximum number of evaluations per variable is set to 50, according to the results by Swarzberg~et~al.~(\citeyear{SW94}). Then, this technique along with its test uses up to $52 \cdot D + 1$ function evaluations.

\vspace{-5pt}
\subsection{The AsBeC algorithm}\label{AsBeC}

The combination of all the previous techniques, integrated into ABC, gives the AsBeC algorithm. Its pseudo-code is presented in Algorithm~\ref{fig.pseudoAsBeC}, indicating the sections in which the new techniques added to ABC are described.
BO and PD are enhancements of ABC regarding bees re-organization, while LI, QP and SGO are hybridizations to estimate the local and global behaviour of the objective function, without changing the mutation equation itself. Thanks to these hybridizations, the bees assume new extra abilities and are called super-bees. This explains the name of the new algorithm: the Artificial super-Bee enhanced Colony (AsBeC).

Actually, the authors experimented also other modifications, e.g. in the mutation equation, and other hybridizations with derivative-free directional search techniques, like the discrete gradient method as in Bagirov et al. (\citeyear{BA08}). In the end, only the variants recognized to have an effective role for the specific goals of this work were implemented within AsBeC.

\begin{algorithm} [!]
	\caption{The AsBeC pseudo-code}
	\label{fig.pseudoAsBeC}
	\begin{algorithmic}[1]
		\LineComment{\{Initialization\}}
		\State{Generate $\max(\textrm{number of seed solutions}, 2 \cdot D +1)$ starting random solutions in the search area (Sec. \ref{Tech4})}
		\State{Evaluate their quality}
		\State{Identify seed solutions as the best starting solutions and assign employees to them}
		\vspace{1 mm}
		\State{Find minimum of the quadratic model near the best seed solution (Sec. \ref{Tech4})}
		\IIf {a new best solution has been found} update best seed solution \EndIIf
		
		\vspace{1 mm}
		\LineComment{\{SGO phase (Sec. \ref{Tech5})\}}
		\IIf {quasi-separability test passed} use interleaved S.T.E.P solver \EndIIf
		\IIf {a new best solution has been found} update best seed solution \EndIIf

		\vspace{1 mm}
		\Repeat
		\LineComment\hspace{12pt}{\{Employees phase\}}
		\ForAll{employees}
		\State Generate new pseudo-random solution near its seed solution
		\State Evaluate the quality of solution
		\IIf{it is better than current employee's seed solution} update seed solution \EndIIf
		\EndFor
		\vspace{1 mm}
		\LineComment{\{Onlookers phase\}}
		\State Assign onlookers to the seed solutions depending on their quality, according to a biased deterministic rule	(Sec. \ref{Tech1})
		\ForAll{onlookers} 
		\ForAll{postponed dance iterations (Sec. \ref{Tech2})
		\State 		Generate new solution near its seed solution, possibly using LI (Sec. \ref{Tech3})
		\State 		Evaluate the quality of solution
		\IIf {it is better than current onlooker's seed solution} update seed solution \EndIIf
		\EndFor}
		\EndFor
		\vspace{1 mm}
		\LineComment{\{Quadratic Prophet phase (Sec. \ref{Tech4})\} }
		\ForAll{seed solution}
		\State 	Collect the samplings nearest to the current seed solution
		\State	Build quadratic model in its neighbourhood and find its minimum
		\IIf {it is better than current seed solution} update seed solution \EndIIf
		\EndFor 
		\vspace{1 mm}

		\vspace{1 mm}
		\LineComment{\{Scout phase\}}
		\IIf {a seed solution is not improved for a limited time} replace it with a new random solution \EndIIf\
		\vspace{1 mm}
		\LineComment{\{Best-so-far\}}
		\IIf {a new best solution has been found} update global best \EndIIf
		\Until{requirements are met}
	\end{algorithmic}
\end{algorithm}

\vspace{-5pt}
\section{Benchmarks}\label{benchmark}

Three different benchmarks for unconstrained problems in the continuous domain are adopted in this paper. The first, Set A, is used to test AsBeC techniques, to provide standard settings for the new algorithm and to compare AsBeC with ABC variants and other notable optimizers. The other two, Set B and C, are popular benchmarks used for competitions in international conferences, herein presented for comparing AsBeC with the most recent mehods presented at the CEC 2015 and GECCO 2015 conferences.

\subsection{Set A}

{\renewcommand\arraystretch{2.0}
	\begin{table} []
		\tiny
		\centering
		\begin{tabular}{|c|l|l|c|c|c|c|}
			\hline
			\textbf{N\#} & \textbf{Name} & \textbf{Formulation} & \textbf{Domain}\ & \textbf{Properties} & \textbf{Minimum} \\
			\hline
			$f_1$ & Sphere & $f(x)=\sum_{i=1}^{D} x_i^2$		&$[-100,100]^D$  &	U, S, Simple bowl & $f(0,...,0)=0$ \\ \hline
			$f_2$ & QuarticR & $f(x)=\sum_{i=1}^{D} i  x_i^4 + rnd[0,1)$ 	& $[-1.28,1.28]^D$ & U, S, Noisy bowl & $f(0,...,0)=0$ \\ \hline
			$f_3$ & Step & $f(x)=\sum_{i=1}^{D}(\lfloor{x_i+0.5}\rfloor)^2$ & $[-100,100]^D$ & U, S, Jagged bowl & $f(0.5,...,0.5)=0$ \\ \hline
			$f_4$ & DixonPrice	& $f(x)=(x_1-1)^2+\sum_{i=2}^{D} i  (2 x_i^2-x_{i-1})^2$ & $[-10,10]^D$ & U, N, Valley & \parbox[c][16pt][c]{1.7cm}{$f( 1,...,2^{-\frac{2^D-2}{2^D}})=0$} \\ \hline
			$f_5$ & Powell	& \parbox[c][22pt][c]{5.3cm}{$f(x)=\sum_{i=1}^{D/4}[(x_{4i-3}+10x_{4i-2})^2+5(x_{4i-1}-x_{4i})^2+(x_{4i-2}-2x_{4i-1})^2+10(x_{4i-3}-x_{4i})^4]$} & $[-4,5]^D$ & U, N, Valley & $f(3,-1,0,1,...,3,-1,0,1)=0$ \\ \hline
			$f_6$ & Rosenbrock	& $f(x)=\sum_{i=1}^{D-1}[100(x_i^2-x_{i+1})^2+(x_i-1)^2]$	& $[-30,30]^D$ & U, N, Valley & $f(1,...,1)=0$ \\ \hline
			$f_7$ & Schwefel1.2	& 	$f(x)=\sum_{i=1}^{D} [\sum_{j=1}^{i} (x_i)]^2$	& $[-100,100]^D$  &  U, N, Ridges & $f(0,...,0)=0$ \\ \hline
			$f_8$ & Schwefel2.22	& 	$f(x)=\sum_{i=1}^{D} |x_i| + \prod_{i=1}^{D} |x_i|$	&  $[-10,10]^D$ &  U, N, Ridges & $f(0,...,0)=0$ \\ \hline
			$f_9$ & Zakharov & $f(x)=\sum_{i=1}^{D} (x_i)^2 + (\sum_{i=1}^{D} 0.5i(x_i))^2+ (\sum_{i=1}^{D} 0.5i(x_i))^4$ & $[-5,10]^D$ & U, N, Large plate & $f(0,...,0)=0$ \\ \hline
			$f_{10}$ & Alpine & 	$f(x)=\sum_{i=1}^{D} | x_i \sin(x_i)+0.1x_i|$	& $[-10,10]^D$ & M, S, Many far minima & $f(0,...,0)=0$ \\ \hline
			$f_{11}$ & Rastrigin &	 $f(x)=10D+\sum_{i=1}^{D} (x_i^2-10\cos(2\pi x_i))$& $[-5.12,5.12]^D$ &  M, S, Several minima & $f(0,...,0)=0$ \\ \hline
			$f_{12}$ & Ackley & $f(x)=-20 \exp(-0.2 \sqrt{\frac{1}{D} \sum_{i=1}^{D} x_i^2}) -\exp(\frac{1}{D} \sum_{i=1}^{D} (\cos(2\pi x_i))) +20 +\exp(1) $	&	$[-32,32]^D$  &  M, N, Narrow hole & $f(0,...,0)=0$ \\ \hline
			$f_{13}$ & Griewank	& $f(x)=\sum_{i=1}^{D} \frac{x_1^2}{4000}-\prod_{i=1}^{D} {\cos{(\frac{x_i}{\sqrt{i}})}}+1$ & $[-600,600]^D$  &  M, N, Several minima & $f(0,...,0)=0$ \\ \hline
			$f_{14}$ & Levy	& \parbox[c][30pt][c]{5.3cm}{$f(x)={\sin^2(\pi\omega_1)}+\sum_{i=1}^{D-1}(\omega_i-1)^2[1+10{\sin^2(\pi \omega_i+1)}]+ (\omega_D-1)^2[1+{\sin^2(2\pi \omega_D)}] $ \newline where $ \omega_i=1+\frac{x_i-1}{4} \forall i$} & $[-10,10]^D$ & M, S, Many far minima & $f(1,...,1)=0$ \\ [2.25ex] \hline
			$f_{15}$ & Penalized	&  \parbox[c][30pt][c]{5.3cm}{$f(x) = \frac{\pi}{D}(10\sin^2(\pi y_1)+(S+(y_D-1)^2))+T$ \newline where $y_i=1+0.25(x_{i}+1)$ , $S=\sum_{i=1}^{D-1} (y_i-1)^2(1+10\sin^2(\pi y_{i+1}))$ and $T=\sum_{i=1}^{D} [(100(x_i-10)^4)(x_i>10)+(100(-x_i-10)^4)(x_i < -10)]$}  &$[-50,50]^D$ & M, N, Many far minima & $f(0,...,0)=0$ \\ \hline
			$f_{16}$ & Penalized2	& \parbox[c][30pt][c]{5.3cm}{$f(x) = 0.1(\sin^2(\pi x_i)+(S+((x(D)-1)^2)(1+\sin^2(2\pi x_d))))+T$ \newline where $S=\sum_{i=1}^{D-1} (x_i-1)^2(1+\sin^2(3 \pi x_{i+1}))$ and $T=\sum_{i=1}^{D} [(100(x_i-5)^4)(x_i>5)+(100(-x_i-5)^4)(x_i < -5)]$} & $[-50,50]^D$ & M, N, Many close minima & $f(0,...,0)=0$ \\ \hline
			$f_{17}$ & Schaffer	& \parbox[c][22pt][c]{5.3cm}{$f(x)=0.5+\frac{(\sin(\sqrt{\sum_{i=1}^{D} x_i^2})-0.5}{(1+0.001\sum_{i=1}^{D} x_i^2)^2}$} &$[-100,100]^D$ & M, N, Few far minima & $f(0,...,0)=0$ \\ \hline
			$f_{18}$ &Whitley	& \parbox[c][22pt][c]{5.3cm}{$f(x)=\sum_{i=1}^{D}\sum_{j=1}^{D}(((100(x_i^2-x_j)^2+(1-x_j)^2)^2)/4000-\cos(100(x_i^2-x_j)^2+(1-x_j)^2)+1)$} & $[-10.24,10.24]^D$ & M, N, Many close minima & $f(1,...,1)=0$ \\
			\hline
		\end{tabular}
		\caption{The analytical Set A. U=unimodal; M=multimodal; S=separable; N=non-separable;}
		\label{Tabfunc}
	\end{table}}

Set A contains unimodal, multimodal, separable and non-separable functions selected among classic analytical problems. Their definitions, search domains and main characteristics are described in Table \ref{Tabfunc}. They are all non-negative and have the global minimum value exactly equal to zero. Functions without an exact analytical characterization of the optimum, e.g., the Styblinski-Tang, are not considered in order to avoid bias on results. Furthermore, all the functions chosen are scalable to any dimension of the variable space. For a rather complete summary of classical analytical functions used for global optimization refer to Jamil and Yang (\citeyear{JA13}). Almost all functions in Table \ref{Tabfunc} are very popular among ABC modifications and usually adopted for this kind of assessment (e.g., Gao and Liu, \citeyear{GL11}).

Set A is heterogeneous and relatively simple, since it does not contain any shifted rotated, highly ill-conditioned or extremely complex function. On the other hand, all the functions comprised present a recognizable overall shape and are potentially solvable using limited FEs. For these reasons, Set A is a good environment to test the techniques integrated into AsBeC.

Each function is investigated with 10 dimensions and a maximum number of FEs equal to $10^3$. Each experiment is repeated 300 times. This settings fit the main target of the AsBeC algorithm, i.e., finding the solution of low-dimensional problems using few FEs. Besides, a long term setting is also introduced, investigating 30 dimensions with maximum number of FEs equal to $5 \cdot 10^4$ and each experiment repeated 30 times. This last setting is in line with other ABC works (e.g., Akay and Karaboga, \citeyear{AK10}) and it is  important to assess premature clustering tendency and refining skills on higher dimensional problems when a lot of FEs are allowed.
 
$10^{-16}$ is considered as the minimum achievable by each function of Set A, i.e., it is the selected \textit{tolerance}, in line with the numerical floating point double precision accuracy of MATLAB. The authors of this paper verified that increasing the value of this tolerance does not affect the qualitative interpretation of results by using the performance metrics defined in Section~\ref{metrics}. Each run of the algorithm over a function stops at the maximum allowed number of FEs. Repeating many times the experiment helps in reducing the influence of the random component.

\subsection{Set B }

The expensive optimization session at CEC conference 2015\footnote{\label{note1}\url{http://www3.ntu.edu.sg/home/EPNSugan/index_files/CEC2015/CEC2015.htm}} is based on a benchmark, here called Set B, developed by Qu~et~al. (\citeyear{QU15}) of 15 functions. The tolerance used is $10^{-8}$ and the number of repetitions is set to $20$. The functions are tested with $10$ dimensions and $500$ $FEs$ and with $30$ dimensions and $1500$ FEs.

Set B has been explicitly targeted to costly optimization and therefore it is appropriate for testing AsBeC. Its functions are shifted rotated, thus analytically not separable, and it includes highly ill-conditioned and very complex problems, also with extremely jagged shapes. Set B is complementary to Set A, since it tests the adaptability of the optimization methods also to these kinds of functions. Notice that the maximum number of admitted FEs in $30$ dimensions is really small ($50 \cdot D$), which is likely to favor a simple exploration of the space.
Further details about the settings and functions can be found in the work by Chen~et~al. (\citeyear{CH15}).

\subsection{Set C}

The noiseless Black-Box Optimization Benchmarking (BBOB)\footnote{\label{note2}\url{http://coco.gforge.inria.fr}} set used at CEC 2015 and at GECCO 2015 conferences, here called Set C, is developed by Hansen~et~al. (\citeyear{HA09a} and \citeyear{HA09b}) and contains 24 functions. The tolerance and the number of repetitions are respectively set to $10^{-8}$ and to $20$. In the present work, the functions are tested with $5$ dimensions and $500$ $FEs$ and with $20$ dimensions and $2000$ $FEs$.

The major advantage of using this benchmark is due to its schematic structure and to the huge amount of literature dealing with it; in fact, at least 151 algorithms have been compared so-far on this benchmark from 2009 to 2015. Set C is not explicitly targeted to costly optimization, since its functions are not meant to be solved within a strictly limited number of FEs. Set C comprises five different groups of functions: (i) separable, (ii) low-conditioned, (iii) high-conditioned, (iv) multi-modal with adequate global structure and (v) multi-modal with weak global structure. This well organized heterogeneous architecture allows to interpret the strengths and weaknesses of each optimization algorithm.
For details about this benchmark refer to Hansen~et~al.~(\citeyear{HA09b}).

\section{Mean logarithmic value metrics}\label{metrics}

Some quantitative metrics are introduced to quickly compare algorithms in terms of the objective value reached after a given amount of FEs. A performance estimator for an algorithm is defined as the Logarithmic Value $LV$:

\begin{equation*}
LV(f,FEs):=\log{\bigg[\frac{(M_{rep}{(f_{best}(f,FEs,rep))-f^*_{best}(f)}}{tolerance }\bigg]}
\end{equation*}

where $rep$ is the run repetition, $FEs$ is the number of performed function evaluations, $f$ is the selected benchmark function, $f^*_{best}$ is the analytical optimal value, $f_{best}$ is the best objective value obtained and $M_{rep}$ represents the median operator over repetitions. The median is chosen in spite of the mean since it is less sensible to outliers. The numerator represents the \textit{residual} and the $LV$ represents the median performance of an algorithm in terms of orders of magnitude missing to the given \textit{tolerance}. A $LV$ is always nonnegative and it reaches zero when the median results of an algorithm are optimal, within the given tolerance.

Three averaged forms of the Logarithmic Value, called Mean Logarithmic Values ($MLVs$), are defined. $MLV_{FEs}$ is an average over the number of benchmark functions $N_f$, $MLV_f$ is an average over the maximum number of function evaluations $FEs_{max}$ and $MLV_A$ is an average over both functions and FEs:

\begin{equation*}
MLV_{FEs}(FEs)=\frac{1}{N_f} \cdot \sum_{f=1}^{N_f} {LV(f,FEs)}
\end{equation*}
\begin{equation*}
MLV_f(f)=\frac{1}{FEs_{max}} \cdot \sum_{FEs=1}^{FEs_{max}} {LV(f,FEs)} \end{equation*}
\begin{equation*} 
MLV_A=\frac{1}{FEs_{max}} \cdot \sum_{FEs=1}^{FEs_{max}} \frac{1}{N_f}  \sum_{f=1}^{N_f} {LV(f,FEs)}
\end{equation*}

$MLV_A$ helps in comparing different algorithms through just one real number. It is meaningful since it considers the performance evolution on the whole set of functions, not only after a specific number of FEs.

The main difference of the introduced metric with respect to a simple mean of the median residuals is the fact that $MLVs$ average the logarithms of the residuals, and not the mean of the actual values. This helps to assess how much the algorithm is getting closer to the optimum in relative terms.

The Logaritmic Value and its averaged versions are here introduced for the first time by the authors of this paper. They will be used for comparisons alongside other classical and benchmark-specific metrics.

\section{Analysis of the improving techniques on Set A}\label{Test}

Having defined the comparison metrics, the main three parameters of the ABC architecture, i.e., the colony size, the number of cycles and the limit parameter, are chosen. The limit parameter is set on $D \cdot SN$ as advised for the original ABC (Karaboga and Akay, \citeyear{KA07b}), while the number of cycles depends on FEs and on the colony size. 
The overall number of agents, $N$, is a key performance-driving parameter for population-based algorithms. Six configurations for the total number of individuals (4, 8, 16, 32, 64 and 128 individuals) are directly investigated and $MLV_A$ is the metric used to choose the best one. For the original ABC the best colony size turns out to be 8 (see Table \ref{tab.N}). This properly tuned ABC is used for the comparisons.

\begin{figure} [!]
    \centering
    \begin{subfigure}[]{0.46\textwidth}
            \includegraphics[width=1\textwidth]{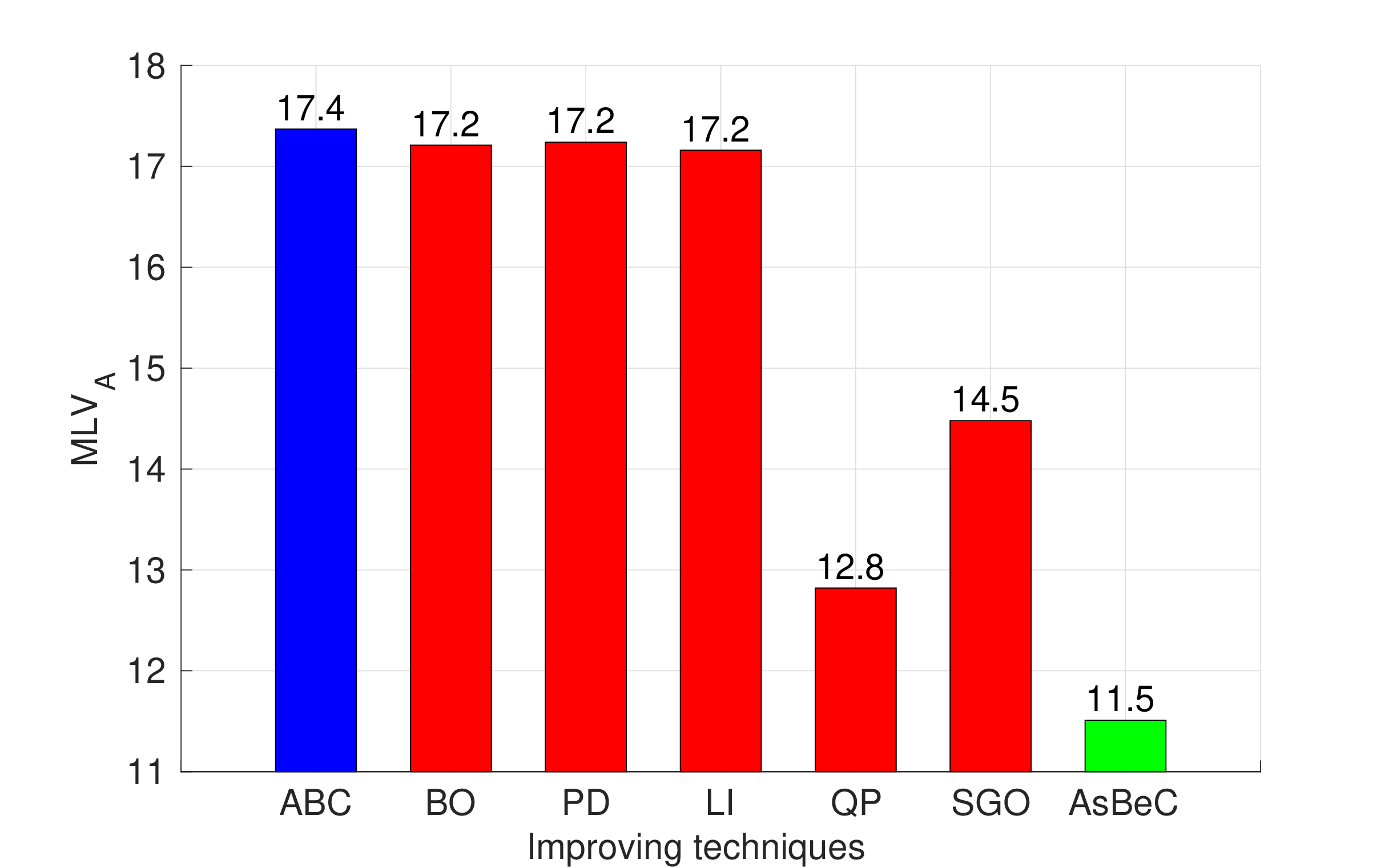}
            \caption{$MLV_A$ for each technique, ABC and AsBeC}
            \label{fig.tech_a}
    \end{subfigure}%
    \quad
    \begin{subfigure}[]{0.49\textwidth}
            \includegraphics[width=1\textwidth]{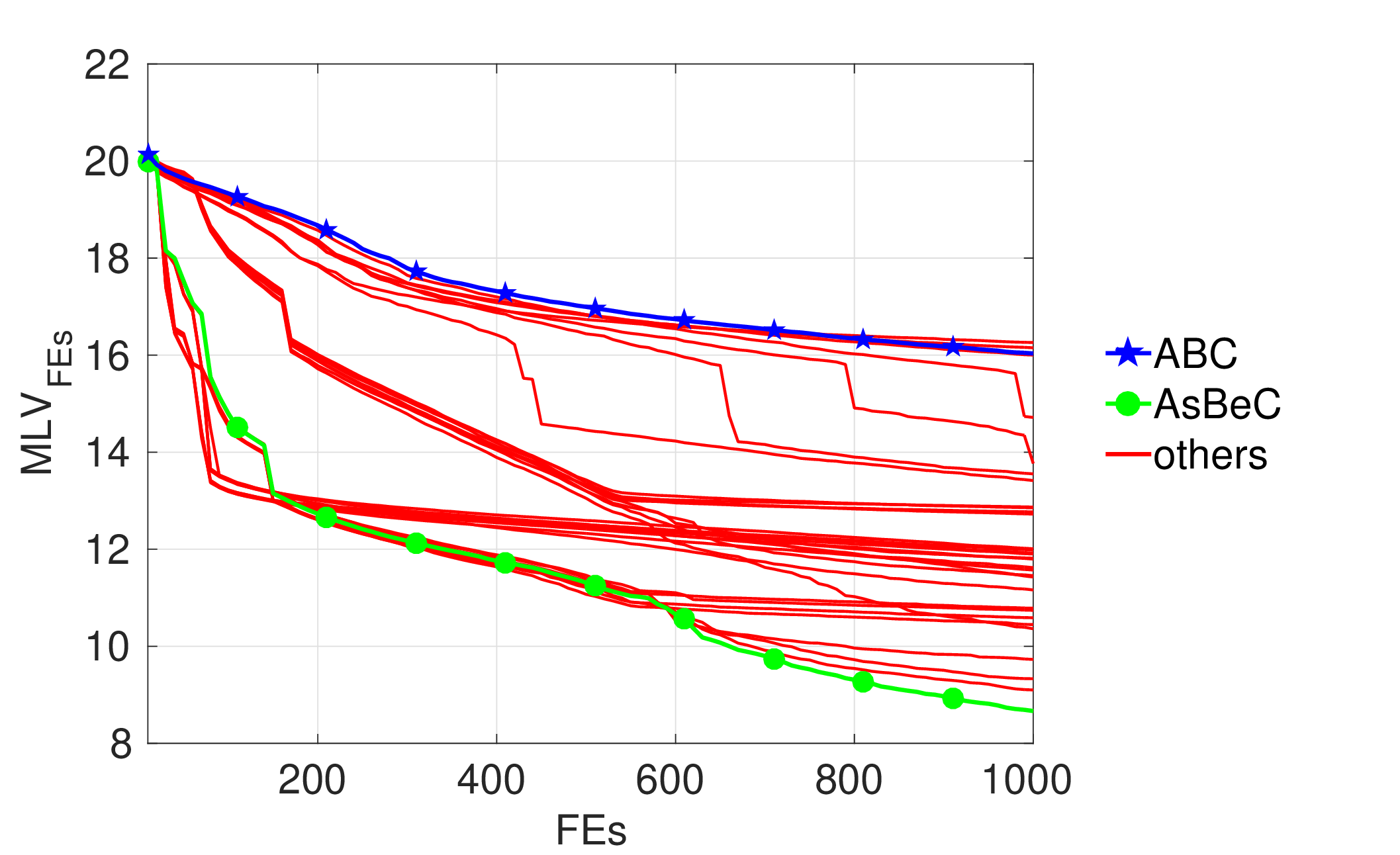}
            \caption{$MLV_{FEs}$ for all the combinations of techniques}
            \label{fig.tech_b}
    \end{subfigure}
	\caption{Validation of the techniques using the $MLV$ metrics}
	\label{fig.tech}
\end{figure}

At this point, the impact of each single technique proposed in Section~\ref{Tech} is evaluated. They are applied one-by-one to the reference ABC configuration. Fig.~\ref{fig.tech_a} illustrates the correspondent $MLV_A$ for each technique and for the basic ABC. All the techniques have positive effects on the algorithm performance. The QP and SGO clearly take a dominant role with respect to the other techniques, that show a more limited impact. 

The differences between the $ABC$ with or without the implemented techniques presented in Fig.~\ref{fig.tech_a} are statistically significant. The $MLV_A$ represents the average of the $LV$ over the evaluations and the functions. A Student significance t-test has been conducted over the two samples (the samples mean will follow a normal distribution by the Central Limit Theorem). The null hypothesis that the $MLV_A$ of the modified algorithm is greater than the one of $ABC$ has been tested against the alternative hypothesis that the mean has positive value. In all the cases, the null hypotheses can be rejected with significance level of 0.05.
 
Then, all possible combinations of techniques are tested in order to capture their mutual interactions. The {$MLV_{FEs}$} comparisons are shown in Fig.~\ref{fig.tech_b}, where each line corresponds to a different combination of techniques. The best configuration (AsBeC line) is the one with lowest $MLV_A$, and it corresponds to the AsBeC ($MLV_A$ of $11.5$, reported in Fig.~\ref{fig.tech_a}). The Quadratic Prophet, based on a trust region approach, is revealed to be impressively suited for this kind of problems, as expected.
 
Once the benefit of the techniques is clear and AsBeC is recognized as the best combination of these techniques, the number of agents is studied for Set A. The best number of bees results 8   for the setting with 10 dimensions and 32 for the setting with 30 dimensions (see Table \ref{tab.N}).
Indeed, other tests on different dimensions showed that the number of individuals should be at least equal to 8 and approximatively proportional to the number of dimensions.
The same AsBeC settings discussed in this section for set A will be also used for set B and C, in order to prove the robustness of this algorithm configuration. 
 
\section{Comparison with other algorithms on Set A}\label{Comparisons}

The AsBeC effectiveness has to be tested when compared to other state-of-the-art methods of the same class suitable for the same goal (see Section ~\ref{Intro}), consistently with what has been done  by Akay and Karaboga (\citeyear{AK10}) or Gao and Liu (\citeyear{GL11}) for their algorithms.
According to the considerations in Section~\ref{ABC},  GABC, JA-ABC5 and RABC are selected as representative of ABC modifications also suited for the specific goal under study.
The authors have selected other three typologies of direct search algorithms for the comparison: FIREFLY (Yang, \citeyear{YA09}), CMA-ES (Hansen, \citeyear{HA06}) and BOBYQA (Powell, \citeyear{PO09}). FIREFLY is a notable competitor of the ABC among the swarm-based techniques, while the evolutionary CMA-ES and the trust region BOBYQA are recognized as state-of-the-art in their fields and have been widely used for engineering applications. 

The implementations used for FIREFLY and CMA-ES are the ones made available by the authors themselves. The MATLAB routine used for BOBYQA is instead the one embedded in the NAG optimization toolbox\footnote{\url{http://www.nag.co.uk/numeric/MB/start.asp}}.

\vspace{-5pt}
\subsection{Choice of parameter values on Set A}\label{Comparisons2}

A key point when comparing different methods consists in tuning the parameters of the optimizers to make them perform at their best.
For all the algorithms considered, the parameter selection is addressed as in Section~\ref{Test} for ABC and AsBeC. Only the number $N$ of individuals for population-based methods and the number $N$ of interpolation samples for BOBYQA are directly investigated on Set A. The range studied for $N$ is a progression of the power of two; its boundaries depend on the specific mechanics driving each algorithm. The best size is selected according to the $MLV_A$ metric and is reported in Table \ref{tab.N}.

All other parameters assume the values suggested in literature as suitable for many problems (Zhu and Kwong, \citeyear{ZK10} for GABC, Kang et al., \citeyear{KAN11} for RABC, Sulaiman et al., \citeyear{SL15} for JA-ABC5, Hansen, \citeyear{HA06} for CMA-ES and Powell, \citeyear{PO09} for BOBYQA). For FIREFLY, the authors use the parameters proposed by Yang (\citeyear{YA10}) for the MATLAB implementation of the algorithm, since they reach better average results on the adopted benchmark with respect to the general suggestions by Yang and He (\citeyear{YH13}). 

\begin{table} [ht!]
	\centering
	\begin{tabular}{|c|c|c|c|c|c|c|c|c|}
		\hline
		Dimensions & ABC & GABC & RABC & JA-ABC5 & FIREFLY & CMA-ES & BOBYQA & AsBeC \\
		\hline
		10D & 8 & 8 & 8 & 16 & 16 & 8 & 32 & 8 \\
		\hline
		30D & 16 & 16 & 16 & 16 & 32 & 32 & 256 & 32 \\
		\hline
	\end{tabular}
	\caption{Best $N$ for each algorithm according to the $MLV_A$ metric }
	\label{tab.N}
\end{table}

\vspace{-30pt}
\subsection{Results of comparison on Set A}\label{Comparisons3}

\begin{figure}[b!]
	\centering
	\begin{subfigure}[]{0.48\textwidth}		
		\includegraphics[trim={8cm 0cm 8cm 1cm},clip,width=\textwidth]{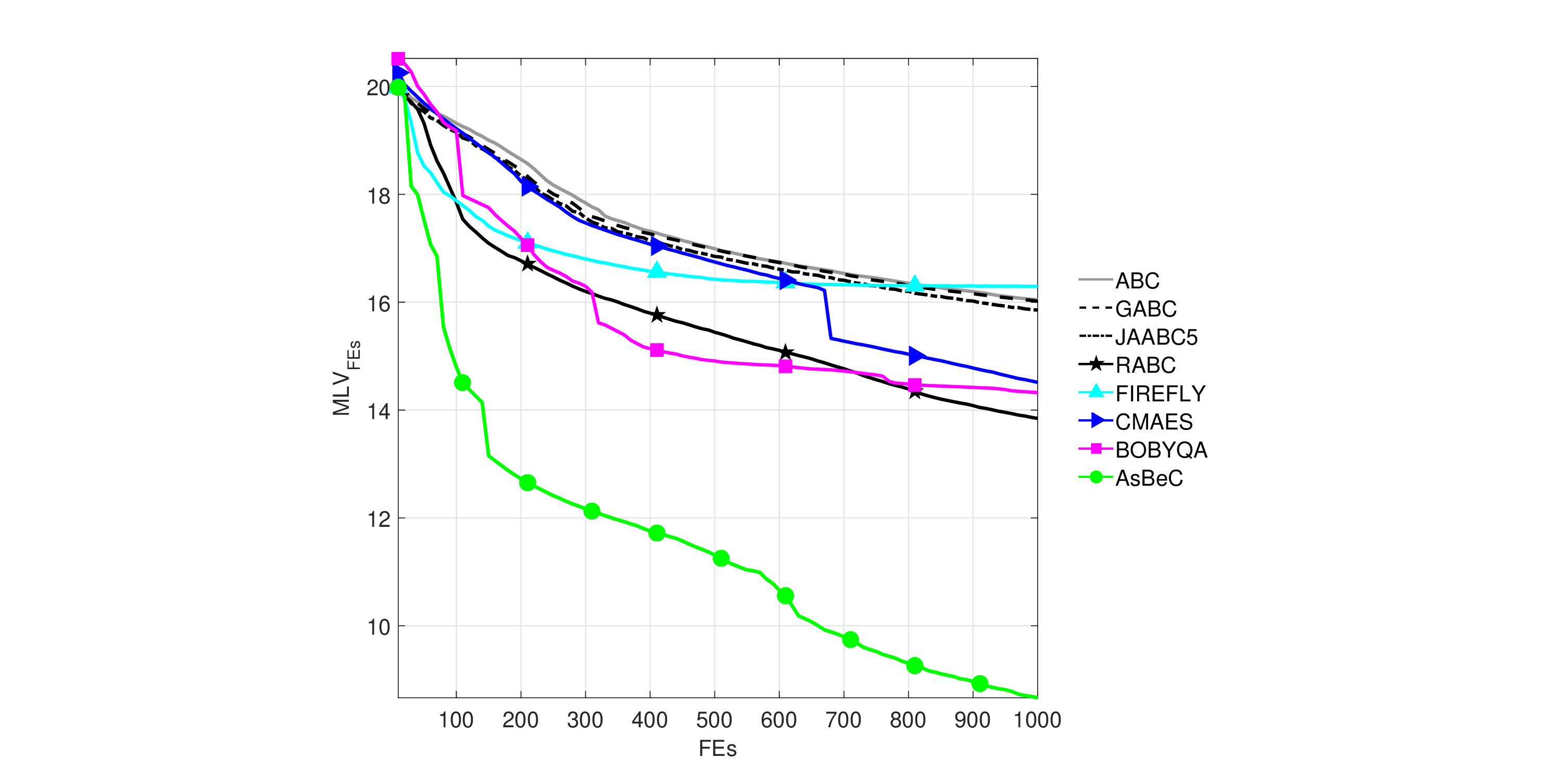}
		\caption{10D}
		\label{fig.setA_a}
	\end{subfigure}
	\begin{subfigure}[]{0.48\textwidth}
		\includegraphics[trim={8cm 0cm 8cm 1cm},clip,width=\textwidth]{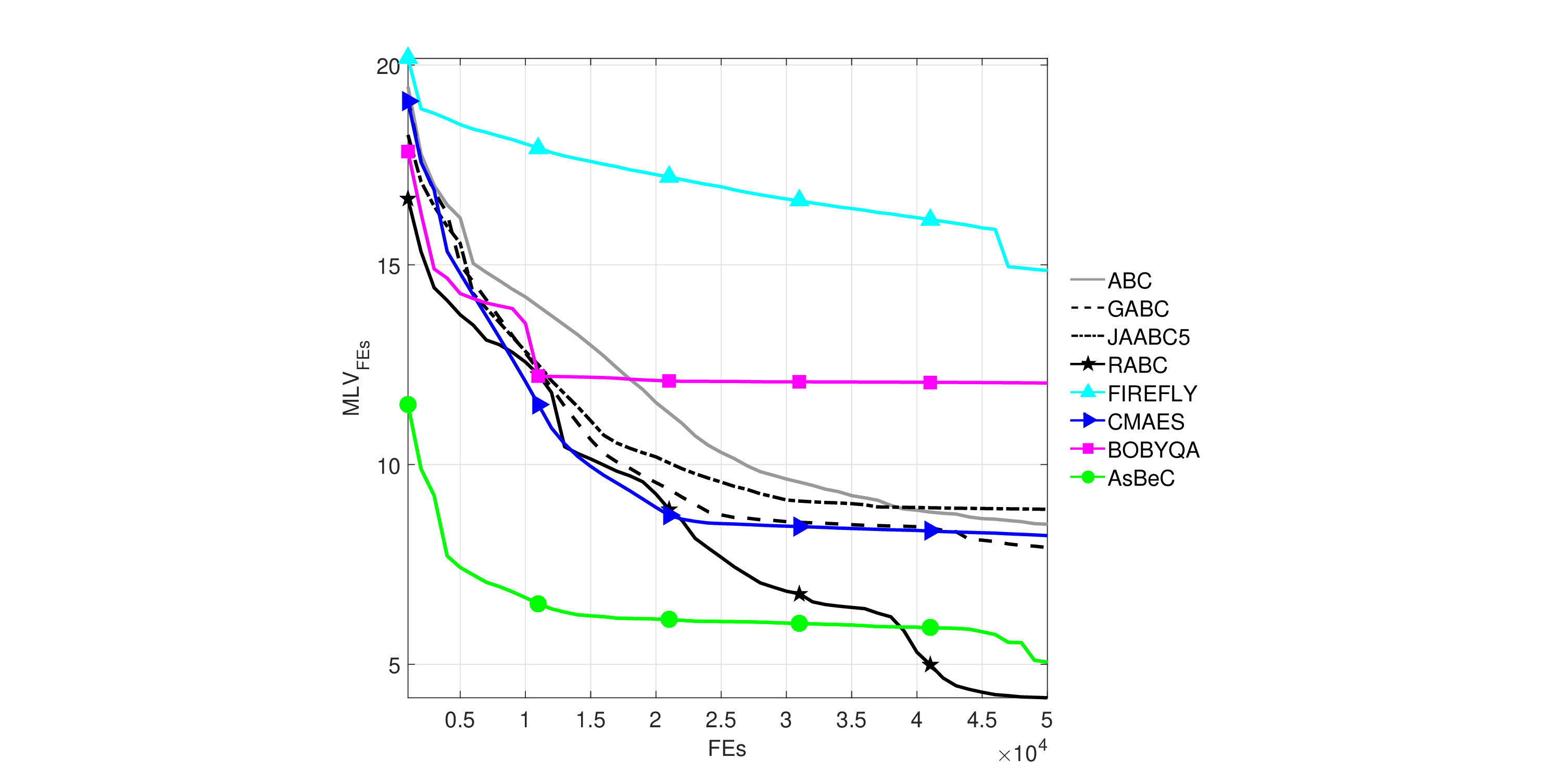}		
		\caption{30D.}
		\label{fig.setA_b}
	\end{subfigure}
	\caption{Set A. Comparison of the algorithms using $MLV_{FEs}$}
	\label{fig.setA}
\end{figure}

The analysis of $MLV_{FEs}$ for 10 dimensions is reported in Figure \ref{fig.setA_a}. It shows that AsBeC is the most effective on average for quickly improving the solution as well as refining it. Referring to function-by-function performance in the electronic appendix, the AsBeC algorithm is able to outperform, or at least to approach, all the others during the entire FEs envelope, on all the benchmark functions. In short, its performance is very good and robust, as it is revealed also by the boxplot comparisons in the electronic appendix.  
None of the algorithms appears remarkably better than AsBeC in any function. AsBeC performs always better than the original ABC and it is the only algorithm able to get close to the global optimum of the Shaffer function. 
Final objective values analysis, reported in the electronic {a}ppendix, recognize AsBeC as the most appropriate for the benchmark. Moreover, following the same procedure adopted in Section~\ref{Test}, the improvements of the $MLV_A$  of AsBeC with respect to the other algorithms show statistical significance.

Considering the other algorithms, the two variants GABC and JA-ABC5 perform slightly better than the original ABC. Instead, RABC is in general faster and more capable also at refining the solution. 
Nonetheless, it is worse than the original ABC in those functions where finding a preferential direction towards the global optimum is particularly difficult. FIREFLY is among the fastest ones in reaching the best regions, but its refining ability is poor. The CMA-ES performance is robust and shows a good growth potential, but it needs more total FEs to improve. BOBYQA presents the most fluctuating results, performing very well on some functions but getting stuck in local minima in many others.

Looking at the overall results, it is clear that they are in general not convergent and in many functions they are still far from the global minimum. Nevertheless, the above comparison is still solid. In fact, in real-like problems where the optimal position is unknown, designers are interested in obtaining the best solution improvement within a limited amount of resources.

Set A with a 30-dimensions setting is not the main focus of the paper, but provides interesting insights on long-term performance. 
The analysis of the MLV in Figure \ref{fig.setA_b} and function-by-function performance reported in the electronic appendix, 
reveal that AsBeC is the most promising for quickly improving the solution also in this kind of problems. In Figure \ref{fig.setA_b}, first data are plotted after $10^3$ FEs to improve visualization, showing that AsBeC reaches a very low value after very few FEs. Moreover, AsBeC is able to reach the given tolerance in 11/18 of the functions and, among the other functions, only Dixon Price and Rosenbrock can be solved by other algorithms.
In short, AsBeC is the best algorithm up to around $4\cdot 10^4$ FEs.
RABC is the most efficient at refining if a large number of FEs is available. 
However, it is clearly much slower than AsBeC. The other two ABC variants and CMA-ES are quite similar and better than ABC, as revealed by the $MLV_{FEs}$. BOBYQA is again among the fastest ones during the very first optimization phases. FIREFLY seems to show a slower convergence speed with respect to the others when dimensionality rises.

\section{Other validations}\label{Validation}

The primary goal of the two validation sets B and C is to assess the robustness of AsBeC, besides its quality, in comparison to the most recent and proficient optimization methods. For this reason, the same AsBeC standard settings discussed in Section~\ref{Test} for Set A will be also used for Set B and C. Following the guidelines of Section \ref{Comparisons2}, the bee colony size is chosen as $N=8$ for $D=5$ and $D=10$,  $N=20$ for $D=20$, and $N=32$ for $D=30$.

\subsection{Results of the comparison on Set B}\label{CEC}

\begin{figure}[t!]
	\centering
	\begin{subfigure}[]{0.48\textwidth}  
		\includegraphics[trim={8cm 0cm 8cm 1cm},clip,width=\textwidth]{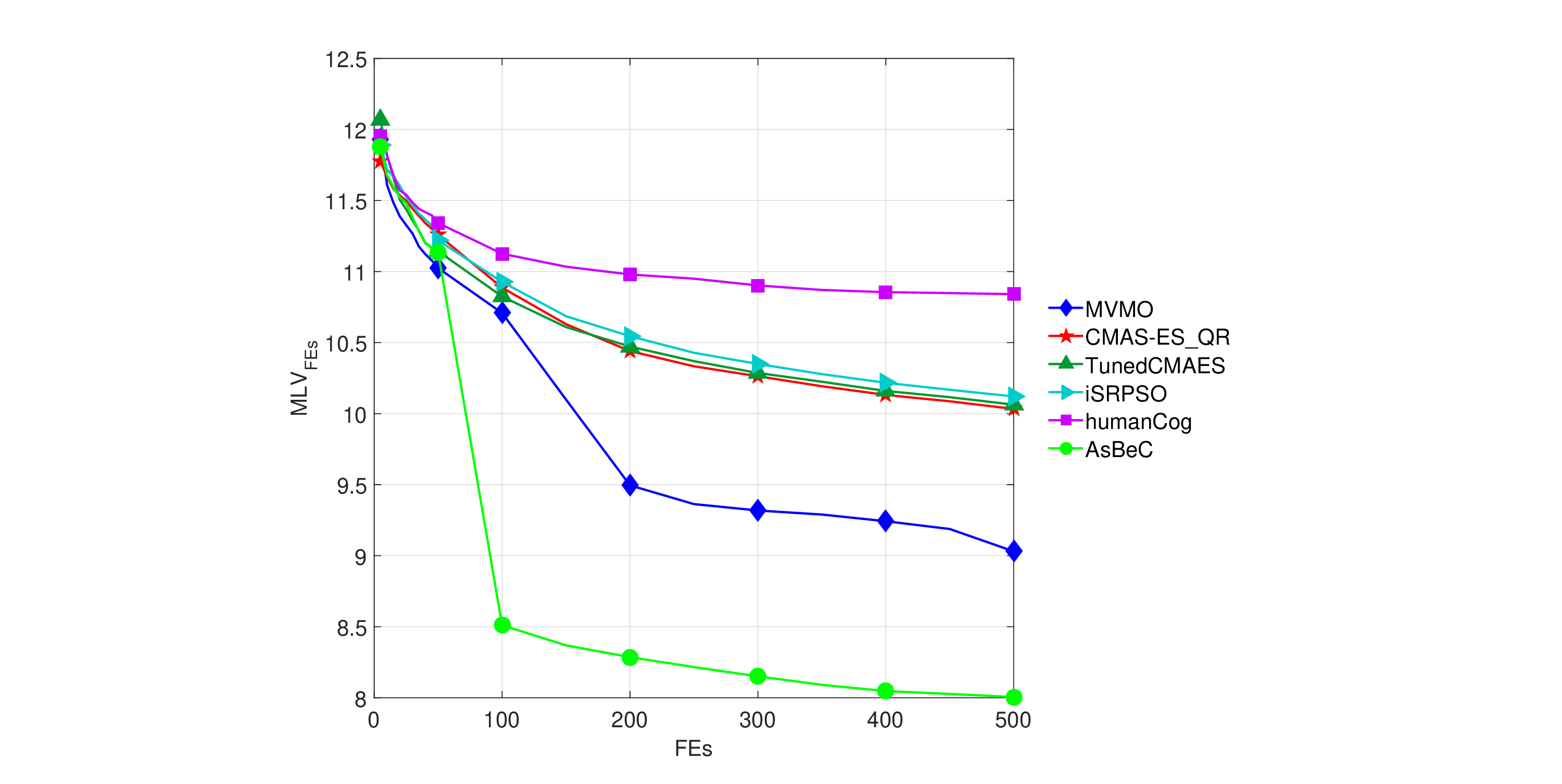}
		\caption{10D}
		\label{fig.setB_a}
	\end{subfigure}
	\begin{subfigure}[]{0.48\textwidth}
		\includegraphics[trim={8cm 0cm 8cm 1cm},clip,width=\textwidth]{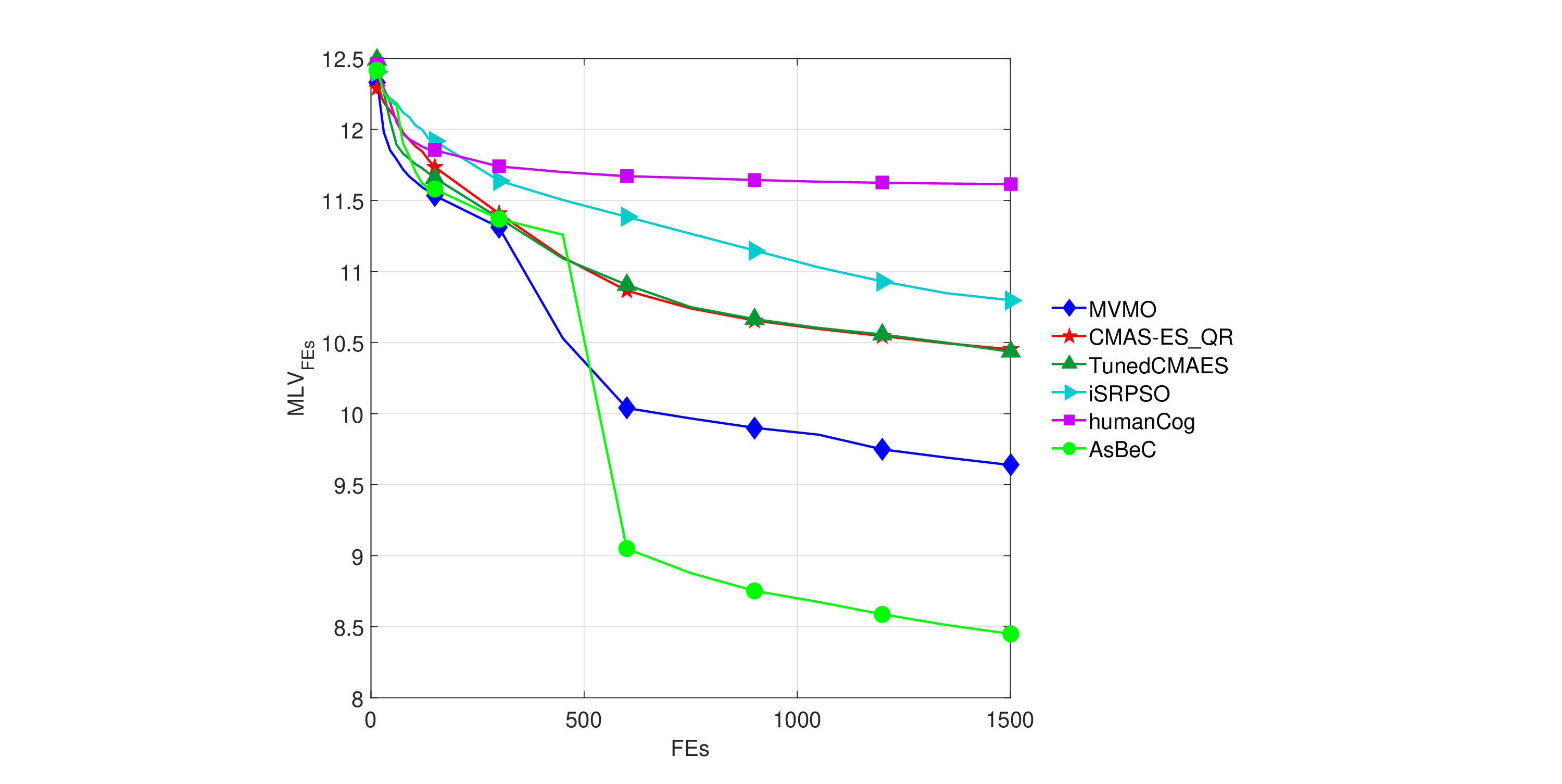}
		\caption{30D}
		\label{fig.setB_b}
	\end{subfigure}
	\caption{Set B. Comparison of the algorithms using $MLV_{FEs}$}
	\label{fig.setB}
\end{figure}

Public results on this benchmark\textsuperscript{\ref{note1}} by Qu~et~al.~(\citeyear{QU15}) include the following algorithms:
\begin{itemize}
	\item MVMO, a population-based stochastic technique with  a mapping function for the offspring;
	\item CMAS-ES QR and TunedCMAES, variants of CMAES for expensive scenarios;
	\item iSRPSO, a PSO implementing a dynamic learning strategy  for velocity updating;
	\item humanCog, a 3-layer architecture that mimics human cognitive behaviour.
\end{itemize}
The $MLV_{FEs}$ comparisons with the AsBeC are presented in Figure \ref{fig.setB}. Extended function-by-function plots of median residual over repetitions as function of FEs and complete tables with mean, median and standard deviation for final achievements are reported in the electronic appendix.
It is evident that AsBeC has a great overall quality, reaching much lower values that the best method, MVMO. AsBeC is able to solve the quadratic problems, even if ill-conditioned, to recognize the quasi-separability of the test function 4, and to approach the best methods in the other functions. Unlike the other algorithms, AsBeC is capable of a stable performance regardless the dimensionality.

Besides, Qu~et~al.~(\citeyear{QU15}) defined an official total scoring for this competition as an average of the median and mean values at the end and at the middle of the computation. The AsBeC total score is one order of magnitude less than MVMO (full table reported in the electronic appendix). Notice that this scoring favors the improvements made on the most difficult function, TF 10. Since both AsBeC and MVMO are able to better minimize this problem, their scores are much lower than those of the other algorithms.

\subsection{Results of comparison on Set C}\label{BBOB}

At the CEC and GECCO conferences in 2015, 26 algorithms have competed on this benchmark\textsuperscript{\ref{note2}}. Indeed, they are tuned variants of few basic methods. Among them, four top categories can be identified:
\begin{itemize}
	\item CMAES derived, including many tuned variants of IPOP-CMAES;
	\item Surrogate based, including the ones that exploit the MATLAB MATSuMoTo Library for metamodels;
	\item DE derived, tuned for cheap, medium, and expensive settings;
	\item Axis-Parallel Brent-S.T.E.P. method, which investigates some variants of the multidimensional Brent-S.T.E.P. method.
\end{itemize}

For each category, the best algorithm in terms of performance and robustness according to the official results of the conference is compared to the standard AsBeC through the $MLV_{FEs}$ metric. These four tuned methods  are the GP5-CMAES, RAND- 2xDefaultMATSuMoTo, R-DE-10e2 and Srr.

\begin{figure}[b!]
	\centering
	\begin{subfigure}[]{0.48\textwidth}
		\label{fig:MLV_5_BBOB}
		\includegraphics[trim={6cm 0cm 8cm 1cm},clip,width=\textwidth]{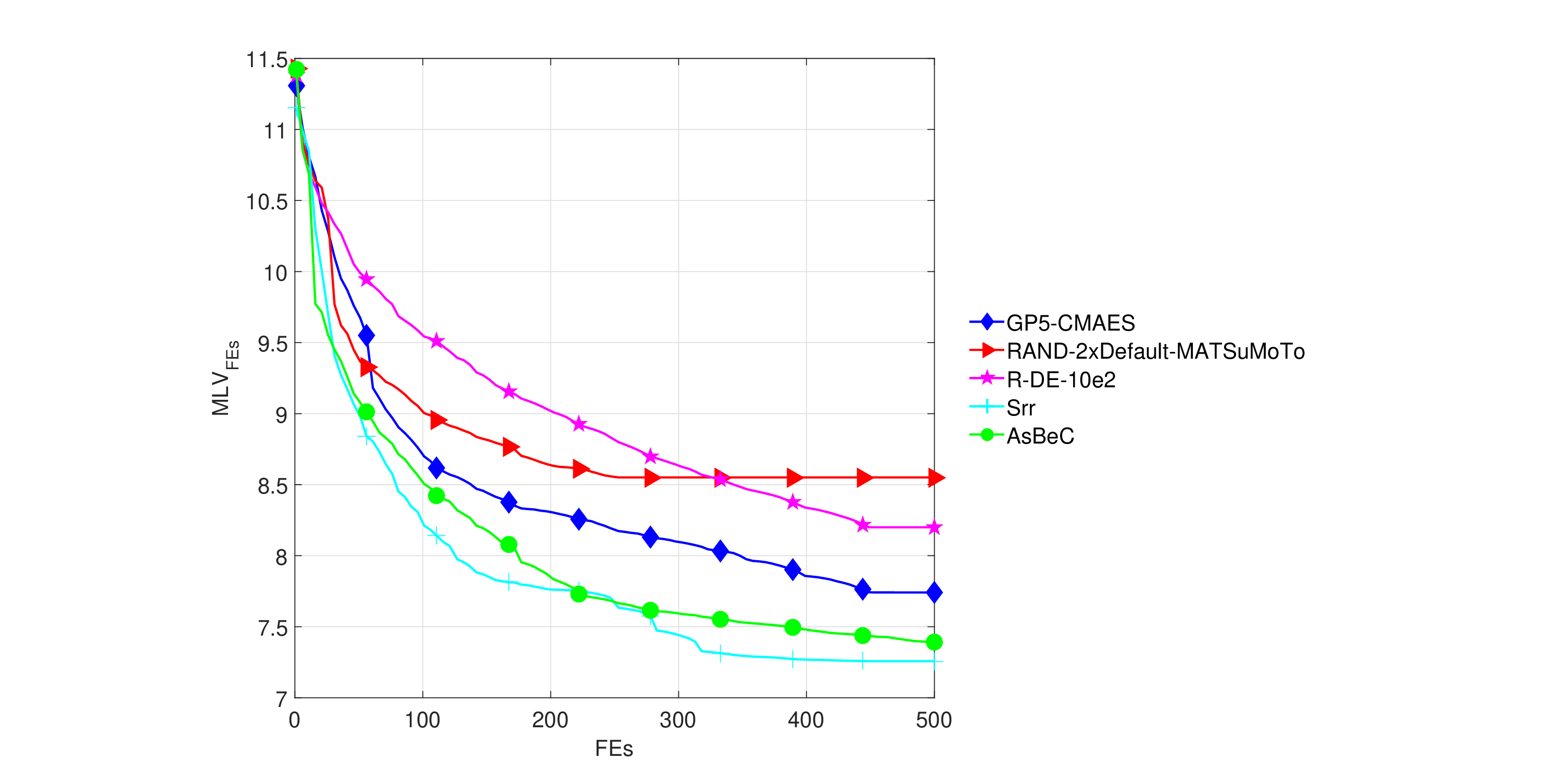}
		\caption{5D}
	\end{subfigure}
	\begin{subfigure}[]{0.48\textwidth}
		\includegraphics[trim={6cm 2cm 8cm 1cm},clip,width=\textwidth]{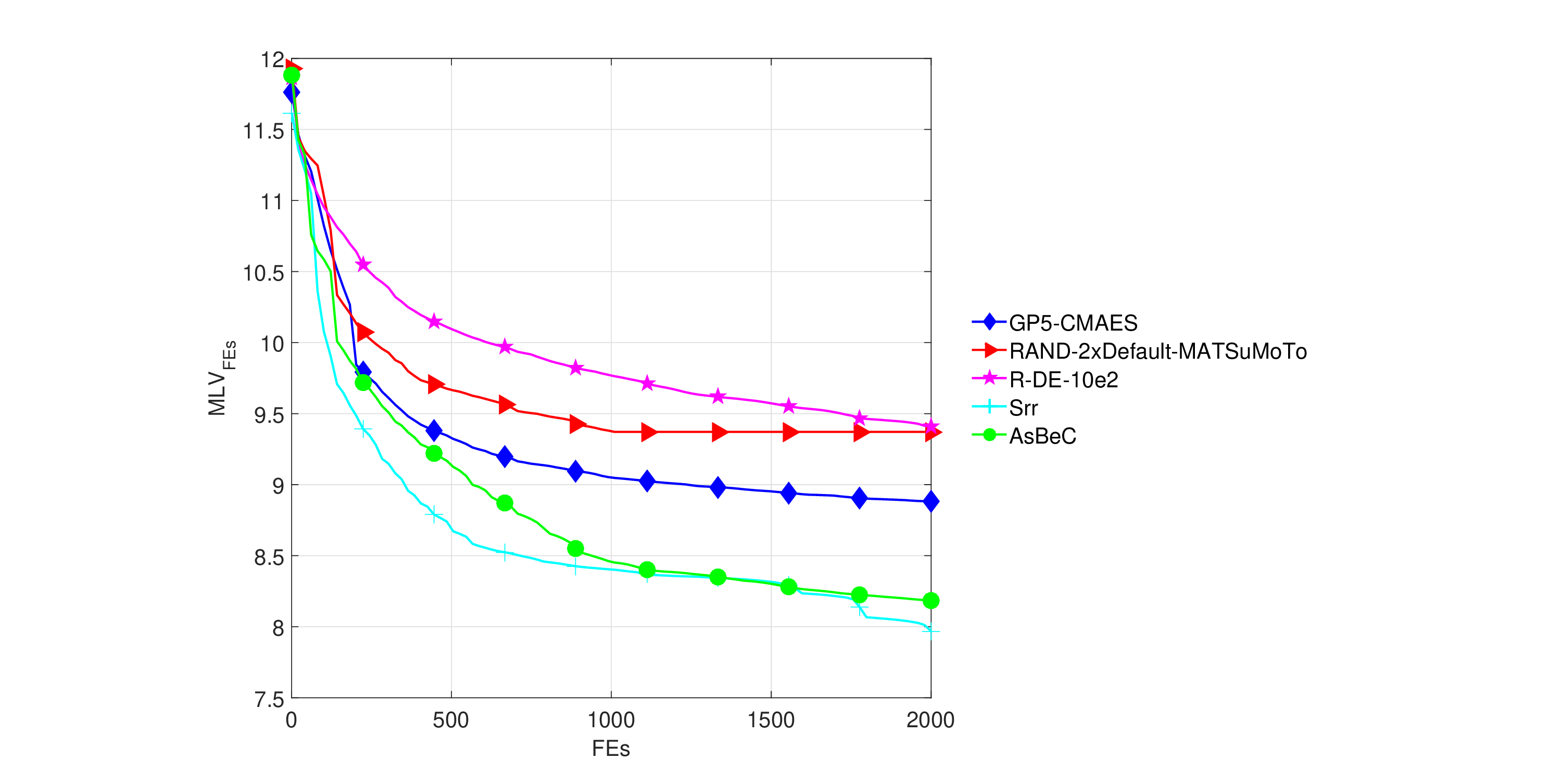}
		\label{fig:MLV_20_BBOB}
		\caption{20D}
	\end{subfigure}
	\caption{Set C. Comparison of the algorithms using $MLV_{FEs}$}
\end{figure}

The $MLV_{FEs}$ metric in Figure~\ref{fig.setB} recognizes the Srr algorithm as the best one, just slightly better than AsBeC. However, function-by-function plots and BBOB official ECDF metric, reported in the electronic appendix, point out that Srr is less robust, especially when addressing non separable multi-modal functions. Srr is the best method for separable problems (f1-5), but it is the worst on some others, reaching premature convergence (e.g., in f7 and f18). Instead, AsBeC and GP5-CMAES offer overall stable high-quality performance. In particular, AsBeC reaches very good results when solving separable, low-conditioned and weak multi-modal problems.

In all the functions, AsBeC is at least close to the best methods, even if no specific tuning has been studied for this set.

\section{Concluding remarks}\label{Conclusion}

Achieving fast and robust improvements in single-objective optimization problems involving expensive analyses means saving precious time and resources. A new swarm-based algorithm hybridized with interpolation strategies, called Artificial super-Bee enhanced Colony (AsBeC) algorithm, is proposed. The new algorithm has been designed for solving expensive problems with low dimensionality, using a limited number of function evaluations. The ambition of this work is to improve the local search, i.e., the exploitation ability, concurrently preserving  the good global attitude of the original ABC, i.e., exploration ability, especially during the first search phases. 

A meaningful metric for comparison is defined, the Mean Logarithmic Value, which takes into account the evolution of the optimization process and the relative distance to the analytical optimum. All the implemented techniques are analysed in order to identify a standard robust setting for them. The standard AsBeC algorithm is compared with ABC, with some of its relevant modifications and with other state-of-the-art direct search algorithms. The same standard AsBeC has also been validated on other benchmarks, and it is compared with some of the latest methods presented at the CEC 2015 and GECCO 2015 conferences, which have specifically been tuned for the benchmarks. The AsBeC algorithm is confirmed to be robust and effective on all the benchmarks.

This promising outcome paves the way for a useful application of the proposed algorithm, especially in engineering. Indeed, the basic principles presented in this paper already showed interesting results in the past when applied by the authors to turbine design (Bertini et al., \citeyear{BE13}). Future works will include comparisons on real-world optimization problems, introduction of parallel strategies and extension to multi-objective problems.

\begin{acknowledgements}
The authors would like to thank GE Avio S.r.l. and its Engineering Technologies department, especially Ing. F. Bertini and Ing. E. Spano. Their collaboration was fundamental for shaping the AsBeC algorithm  within an industrial application framework. We are also grateful to Prof. E. Benini (Universit\`a di Padova) for his comments and reviews.
\end{acknowledgements}

\bibliographystyle{plainnat}

\end{document}